\documentclass[a4paper,12pt]{article}
    \usepackage[top=2.5cm,bottom=2.5cm,left=2.5cm,right=2.5cm]{geometry}
    \usepackage{cite, amsmath, amssymb}
\usepackage{tikz}
\usepackage{amsmath}

\begin{document}
\begin{center}
{\LARGE\bf More on Equienergetic Threshold graphs}
\end{center}
\begin{center}
{\large \bf Fernando Tura}
\end{center}
\begin{center}
\it Departamento de Matem\'atica, UFSM, Santa Maria, RS, 97105-900, Brazil
\end{center}
\begin{center}
{\tt ftura@smail.ufsm.br}
\end{center}



\newtheorem{Thr}{Theorem}
\newtheorem{Pro}{Proposition}
\newtheorem{Que}{Question}
\newtheorem{Con}{Conjecture}
\newtheorem{Cor}{Corollary}
\newtheorem{Lem}{Lemma}
\newtheorem{Fac}{Fact}
\newtheorem{Ex}{Example}
\newtheorem{Def}{Definition}
\newtheorem{Prop}{Proposition}
\def\floor#1{\left\lfloor{#1}\right\rfloor}

\newenvironment{my_enumerate}{
\begin{enumerate}
  \setlength{\baselineskip}{14pt}
  \setlength{\parskip}{0pt}
  \setlength{\parsep}{0pt}}{\end{enumerate}
}

\newenvironment{my_description}{
\begin{description}
  \setlength{\baselineskip}{14pt}
  \setlength{\parskip}{0pt}
  \setlength{\parsep}{0pt}}{\end{description}
}


\begin{abstract} The energy of a graph is defined  as the sum the absolute values of  the eigenvalues of its adjacency matrix. A threshold graph $G$ on $n$ vertices is coded  by a binary sequence of length $n.$ In this paper we answer a question posed by Jacobs et al. [Eigenvalues and energy in threshold graphs, Linear Algebra Appl. 465 (2015) 412-425], by giving an infinite sequences of threshold graphs having the same energy but differ to the complete graph's energy.
\end{abstract}

\baselineskip=0.30in


\section{Introduction}
\label{intro}

Let  $G= (V,E)$ be an undirected graph with vertex set $V$ and edge set $E,$ without loops or multiple edges. 
The {\em adjacency matrix} of $G$, denoted by $A=[a_{ij}]$, is a matrix whose rows and columns are indexed by the vertices of $G$, and is defined to have entries $a_{ij}=1$ if and only if $v_i$ is adjacent to $v_j,$ and $a_{ij}=0$ otherwise.

If $G$ is a graph of order $n,$ its {\em energy} is defined as
\begin{equation}
\label{eq1}
 E(G)= \sum_{i=1}^n |  \lambda_i|
 \end{equation}
 where $\lambda_i$ are the eigenvalues of its adjacency matrix.
There are many results on energy and its applications in several areas, including in chemistral
 see  \cite{Gutman2012} for more details and the references \cite{Li, Li2, Li3, Li4, Li5, Li6}.

It is well known that the complete graph $K_n$ has $E(K_n) = 2n-2.$ According \cite{Gutman2015}, a graph $G$ on $n$ vertices is said to be  {\em borderenergetic} if its energy satisfies
$E(G)= 2n-2.$ In \cite{JTT2015} considered the eigenvalues and energies of threshold graphs. For each $n\geq 3,$ they determined  $n-1$ threshold graphs on $n^2$ vertices, pairwise non-cospectral and equienergetic to the complete graph $K_{n^2}.$ Recently, \cite{Hou} generalized the results presented in \cite{JTT2015}.

Finding all equienergetic threshold graphs was the motivation of this work.  
Even though we have not succeeded, we did answered the question posed by Jacobs et al. \cite{JTT2015}: if there are threshold graphs having the same energy but differ to the complete graph's energy. The answer is affirmative. Indeed, we exhibhit infinite sequences of noncospectral  and equienergetic threshold graphs.


The paper is organized as follows.
In Section \ref{Sec2}  we show the representation of a threshold graph by binary sequence and some known results. In Section \ref{Sec3} we present an infinite sequences
of equinegetic threshold graphs having the same energy but differ to the complete graph's energy.

\section{Preliminaries}
\label{Sec2}
We recall the definition of threshold graphs from \cite{JTT2015}.   
A threshold graph $G$ on $n$ vertices is coded by a binary sequence $(b_1 b_2 \ldots b_n).$
Here  $b_i=0$ if an isolated  vertex $v_i$ is added and $b_i=1$  if $v_i$ was added as dominating vertex. The choice of digit associated to $v_1$ is arbitrary and we use it as $b_1=0.$

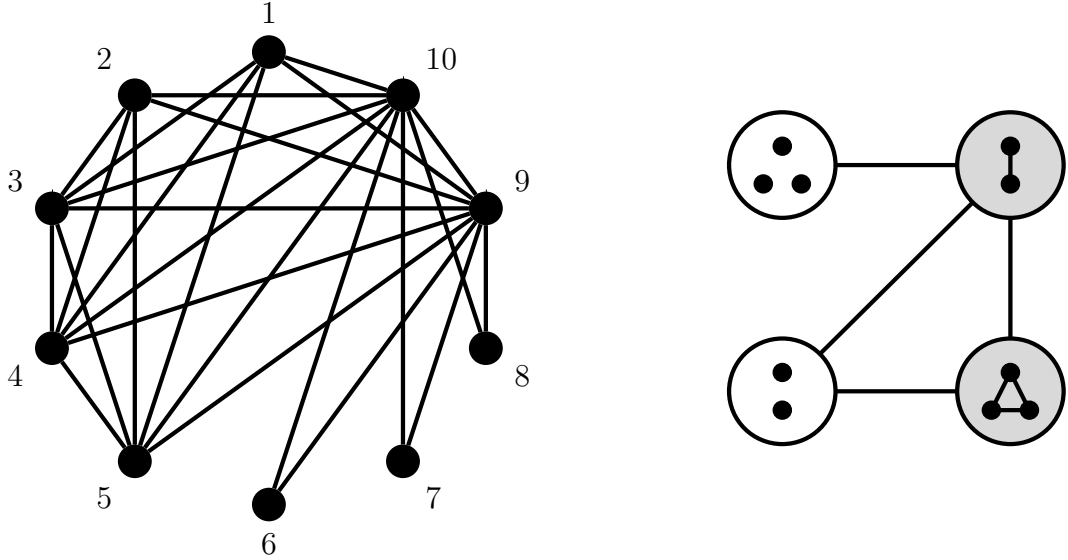
\begin{figure}[h!]
       \begin{minipage}[c]{0.25 \linewidth}
\begin{tikzpicture}[ultra thick]
  [scale=0.5,auto=left,every node/.style={circle}]
  \foreach \i/\w in {1/,2/,3/,4/,5/,6/,7/,8/,9/,10/}{
    \node[draw,circle,fill=black,label={360/10 * (\i - 1)+90}:\i] (\i) at ({360/10 * (\i - 1)+90}:3) {\w};} 
  \foreach \from in {3,4,5,9,10}{
    \foreach \to in {1,2,...,\from}
      \draw (\from) -- (\to);}
\end{tikzpicture}
       \end{minipage}\hfill
       \begin{minipage}[l]{0.4 \linewidth}
\begin{tikzpicture}[ultra thick]
\draw[fill=white ](0,0) circle [radius=0.7];
\draw[fill](0, 0.25) circle[radius =0.1];
\draw[fill](0, -0.25) circle[radius =0.1];

\draw[fill=white ](0,3) circle [radius=0.7];
\draw[fill](0, 3.25) circle[radius =0.1];
\draw[fill](-0.25, 2.75) circle[radius =0.1];
\draw[fill](0.25, 2.75) circle[radius =0.1];

\draw(0.7, 3)--(2.3,3);
\draw(3, 2.3)--(3,0.7);
\draw(0.7,0)--(2.3,0);

\draw(0.5, 0.5)--(2.5,2.5);


\draw[fill=gray!30! white](3,3) circle[radius=0.7];
\draw(3, 3.25)--(3,2.75);
\draw[fill](3, 3.25) circle[radius =0.1];
\draw[fill](3, 2.75) circle[radius =0.1];

\draw[fill=gray!30!white ](3,0) circle [radius=0.7];
\draw[fill](3, 0.25) circle[radius =0.1];
\draw[fill](2.75, -0.25) circle[radius =0.1];
\draw[fill](3.25, -0.25) circle[radius =0.1];
\draw(3, 0.25)--(2.75,-0.25);
\draw(3, 0.25)--(3.25,-0.25);
\draw(2.75, -0.25)--(3.25,-0.25);

\end{tikzpicture}
       \end{minipage}
       \caption{A threshold graph and its partitioned representation}
       \label{fig1}
\end{figure}

During this paper, we denote $G= (0^{a_1} 1^{a_2} \ldots 0^{a_{n-1}} 1^{a_n})$  a connected threshold graph $G$ 
where each $a_i$ is a positive integer.  As illustration  the Figure \ref{fig1} shows a threshold graph $G=(0^2 1^3 0^3 1^2)$ and its partitioned representation.


An interesting result about the multiplicity of eigenvalues $0$ and $-1$ of a threshold  graph $G$ was presented in \cite{JTT2015}.
Let $m_{0}(G)$ and $m_{-1}(G)$ denote the multiplicity of the eigenvalue $0$ and $-1,$ respectively,  they can be obtained  directly from its binary sequence, that is
\begin{Thr}
\label{main1} For a connected threshold graph  $G= (0^{a_1} 1^{a_2} \ldots 0^{a_{n-1}} 1^{a_n})$
where each $a_i$ is a positive integer. Then   
\begin{itemize}
 \item[i.]  $m_{0}(G) = \sum_{i=1}^{\frac{n}{2}} (a_{2i-1} -1).$ 

\item[ii.] $m_{-1}(G) = \left\{
\begin{array}{lr}
\sum_{i=1}^{\frac{n}{2}} (a_{2i} -1) & \mbox{if $a_1 > 1$   }\\
 1+ \sum_{i=1}^{\frac{n}{2}} (a_{2i} -1) & \mbox{if $a_1=1.$  }
\end{array} \right.$
\end{itemize}
\end{Thr}

Let $[n]=\{1,2,\ldots n\}$, and let $I_{n,l}$ be the set of increasing sequences of length $l,$ alternating even and odds numbers such that the last term has the same than parity $n$. For instance
$$I_{7,4}=\{(2,3,4,5),(2,3,4,7),(2,3,6,7), (2,5,6,7),(4,5,6,7) \},$$ while 
$$I_{6,4}=\{(1,2,3,4),(1,2,3,6),(1,2,5,6),(1,4,5,6),(3,4,5,6) \}.$$
Given a sequence $\mathbf{t}=(t_1,t_2,\ldots,t_l)$ we denote  $a_\mathbf{t}=a_{t_1} a_{t_2}\cdots a_{t_l}$. Based in this notation,
if $(a_1, a_2, \ldots, a_n)$  is a fixed sequence of positive integers, we define the following parameter 

  $$   \gamma _{n}(l)=  \left\{\begin{array}{ccc}
           \sum\limits_{\mathbf{t}\in I_{n,l}}  a_\mathbf{t}& if &1\leq l \leq n \\
                1  & if &    l=0      \\
   \end{array}\right.$$



The following theorem gives an explicit formula for computing the characteristic polynomial of a threshold graph
from its binary sequence.  This result will play an important role in the sequel.

\begin{Thr}[\cite{JOF} Theorem 5]
\label{main2}
Let $G= (0^{a_1} 1^{a_2} \ldots 0^{a_{n-1}} 1^{a_n})$ be a connected threshold graph where each $a_i$ is a positive integer.
Let $m_0(G)$ and $m_{-1}(G)$ be the multiplicities of eigenvalues $0$ and $-1$ of $G,$ respectively. The characteristic polynomial of $G,$ denoted by $P_G(x),$ is  
$$P_G(x) =(-1)^{\sum a_i} x^{m_0(G)} (x+1)^{m_{-1}(G)} Q_n(x), \hspace{0,5cm} where $$ 
$$Q_n(x)=x^{r_0} \sum_{k=0}^{m} (-1)^{m-k}x^k y^k \gamma_n (n-2 k-r_0)  + x^{r_1} \sum_{k=0}^{m-r_1} (-1)^{m-k}x^k y^k \gamma_n(n-2k-r_1),$$
with $r_0, r_1\in\{0,1\},$ such that, $n=2m+r_0, r_1\equiv r_0+1 (mod \hspace{0,2cm}2)$ and $y=x+1.$
\end{Thr}

\begin{Ex}
We apply the formula given in Theorem \ref{main2} to the threshold graph $G=(0^{a_1} 1^{a_2} 0^{a_3} 1^{a_4})$ with $a_1>1.$ According Theorem \ref{main1} the multiplicities of $0$ and $-1$ are
$m_0(G) = \sum_{i=1}^{2} (a_{2i-1} -1)$ and $m_{-1}(G) = \sum_{i=1}^{2} (a_{2i} -1).$ The rest of eigenvalues are the roots of the polynomial $Q_4(x),$ where
\begin{align*}
Q_4(x) &=x^{0} \sum\limits_{k=0}^{2} (-1)^{1-k}x^k y^k \gamma_4 \text{\footnotesize $(4-2 k-0)$}  +\,  x^{1} \sum\limits_{k=0}^{2-1} (-1)^{1-k}x^k y^k \gamma_4\text{\footnotesize $(4-2k-1)$}\\
&=  x^0 \left( -\gamma_4 (4) +x y \gamma_4 (2) - x^2 y^2 \gamma_4(0)  \right)  +  x^{1} \big(-\gamma_4(3)+x y \gamma_4(1)\big)\\
& =x^2y^2  -(a_2 +a_4)x^2y +(a_1a_2 +  a_1 a_4+a_3a_4)xy +(a_2 a_3 a_4)x  -a_1a_2a_3 a_4. \\
\end{align*}
Therefore the $P_G(x)$ of $G=(0^{a_1} 1^{a_2} 0^{a_3} 1^{a_4})$ is given by
\begin{align*}
P_G(x) &= (-1)^{\sum a_i} x^{m_0(G)} (x+1)^{m_{-1}(G)}   \{x^2 (x+1)^2 -x^2(x+1)(a_2+a_4)    \\
&-(a_1a_2 +  a_1 a_4+a_3a_4)x(x+1) +(a_2 a_3 a_4)x  -a_1a_2a_3 a_4 \}.\\
\end{align*}
\end{Ex}

\section{Main results}
\label{Sec3}

\begin{Lem}
\label{lem1}
For positive integer $i,$ the characteristic polynomial of threshold graph $G=(0^{2i+1} 1^{3i+3} 0^{2i+1} 1^{2i})$ is, to within a sign,
\begin{equation}
\label{eq2}
P_G(x) = x^{4i} (x+1)^{5i+1} (x+2i+1)(x^3  -(7i+2)x^2 -(7i+3)x +12i^3 +18i^2 +6i)
 \end{equation}
\end{Lem}
\noindent{\bf Proof:}  Let $G=(0^{2i+1} 1^{3i+3} 0^{2i+1} 1^{2i})$ be a threshold graph. For computing the characteristic polynomial of $G,$
we use the formula given in the Example 1, with  $a_1=2i+1, a_2= 3i+3, a_3= 2i+1$ and $a_4=2i.$  By Theorem \ref{main1}, we have
$m_0(G)= a_1 + a_3 -2 = 2i+1 +2i+1 -2= 4i$ and $m_{-1}(G) = a_2 + a_4 -2 = 3i +3+2i -2 = 5i+1.$
According Theorem \ref{main2}, the $P_G(x) =x^{4i} (x+1)^{5i+1} Q_4(x)$ where
$Q_4(x) = x^2 (x+1)^2 -x^2(x+1)(5i+3)-(14i^2 +13i +3)x(x+1) +( 12i^3 +18i^2+6i)x -(24i^4 +48i^3+30i^2+6i).$
Factoring the polynomial $Q_4(x),$ the result follows. $\hspace{0,5cm} \square$

\begin{Lem}
\label{lem2}
For positive integer $i,$ the characteristic polynomial of threshold graph $G^{\prime}=(0^{2i+2} 1^{3i} 0^{2i+1} 1^{2i+2})$ is, to within a sign,
\begin{equation}
\label{eq3}
P_{G^{\prime}}(x) = x^{4i+1} (x+1)^{5i} (x+2i+2)(x^3  -(7i+2)x^2 -(7i+3)x +12i^3 +18i^2 +6i)
 \end{equation}
\end{Lem}
\noindent{\bf Proof:}  The proof follows similarly from the Lemma \ref{lem1}.


\begin{Thr}
\label{main3}
The $n$-vertex graphs  $G=(0^{2i+1} 1^{3i+3} 0^{2i+1} 1^{2i})$ and  $G^{\prime}=(0^{2i+2} 1^{3i} 0^{2i+1} 1^{2i+2})$ are noncospectral and equienergetic graphs  if $n=9i+5, i\geq 1.$
\end{Thr}
\noindent{\bf Proof:} It is easy to see from equations (\ref{eq2}) and (\ref{eq3}) that $-2i-1$ and $-2i-2$ are eigenvalues of $G$ and $G^{\prime},$ respectively. Since that $G$ and $G^{\prime}$  share the remainder eigenvalues, this shows that both graphs are noncospectral. Now we will show they have the same energy.  Let $spect(G) = \{ 0, -1, -2i-1, \lambda_1, \lambda_2, \lambda_3 \}$ and $spect(G^{\prime}) = \{0, -1, -2i -2, \lambda'_1, \lambda'_2, \lambda'_3\}$ be the spectrum of $G$ and $G^{\prime}$ without multiplicity included, respectively.  From equation (\ref{eq1}), we have
\begin{equation}
\label{eq4}
E(G) = 5i+1 + 2i+1 +\sum_{i=1}^3 |\lambda_i|
\end{equation}
and
\begin{equation}
\label{eq5}
E(G^{\prime}) = 5i+ 2i+2 +\sum_{i=1}^3 |\lambda'_i|
\end{equation}
Since that $\lambda_i = \lambda'_i$ for $i \in \{1,2,3\},$ from (\ref{eq4}) and (\ref{eq5}) follows that $G$ and $G^{\prime}$ are equienergetic graphs.$\hspace{0,3cm} \square$

\begin{Lem}
\label{lem3}
Let $i\geq 1$ and let $\lambda_1 \leq \lambda_2 \leq \lambda_3$ be the roots of the cubic polynomial
\begin{equation}
\label{eq6}
q(x)=  x^3  -(7i+2)x^2 -(7i+3)x +12i^3 +18i^2 +6i
\end{equation} 
Then $-2i-1 < \lambda_1 <0$ and $\lambda_2 >0.$
\end{Lem}
\noindent{\bf Proof:} We note that $q(0)= 12i^3+18i^2+6i >0.$ Since $q(x)$ is a cubic polynomial, it suffices to verify that $q(-2i-1) < 0.$
By direct calculus, we have that $q(-2i-1) = -24i^3-16i^2+6i < 0,$ for $i\geq 1,$ and so $-2i-1 < \lambda_1 < 0.$  Similarly, we prove that $\lambda_2 >0,$
and the result follows. $\square$

\begin{Cor}
\label{cor1}
The $n$-vertex graphs  $G=(0^{2i+1} 1^{3i+3} 0^{2i+1} 1^{2i})$ and  $G^{\prime}=(0^{2i+2} 1^{3i} 0^{2i+1} 1^{2i+2})$ have the same energy but differ to the complete graph's energy. 
Furthermore, $E(G) = E(G^{\prime}) < E(K_n) = 18i +8,$  if $n=9i+5, i\geq 1.$
\end{Cor}
\noindent{\bf Proof:}  From Theorem \ref{main3}, we have that $G=(0^{2i+1} 1^{3i+3} 0^{2i+1} 1^{2i})$ and  $G^{\prime}=(0^{2i+2} 1^{3i} 0^{2i+1} 1^{2i+2})$ are noncospectral
and equienergetic graphs. Now we will show the complete graph's energy is a bound for the energy of them. Let $\lambda_1 \leq \lambda_2 \leq \lambda_3$ be the roots of the cubic equation (\ref{eq6}).
From Lemma \ref{lem3} it follows that $-2i-1 < \lambda_1$ and $0< \lambda_2.$ Since $\lambda_1 + \lambda_2 + \lambda_3 = 7i+2,$ we have $|\lambda_1| + |\lambda_2| + |\lambda_3|= -\lambda_1 + \lambda_2 + \lambda_3 = \lambda_1 + \lambda_2 +\lambda -2\lambda_1 \leq 7i+2 +2(2i+1) =11i+4.$ Replacing this bound in the equation (\ref{eq4}), we have
$$ E(G) = 7i +2 + \sum_{i=1}^3 |\lambda_i| \leq 7i+2 + 11i+4 = 18i +6 < 18i+8 = E(K_n)$$
and hence the result follows. $\square$

The proof of the following results are similar to others above, then we will omite them.

\begin{Lem}
\label{lem4}
For positive integer $i,$ the characteristic polynomial of threshold graph $G=(0 1^{2i+1}0^i 1^{2i+2} 0^{2i+1} 1^{2i})$ is, to within a sign,
\begin{equation}
\label{eq7}
P_G(x) = x^{3i-1} (x+1)^{6i+1} (x+2i+1)(x^4 -(8i+2)x^3 -(-8i^2+4i+3)x^2 -(-8i^3 -20i^2-8i)x -8i^4-12i^3-4i^2)
 \end{equation}
\end{Lem}

\begin{Lem}
\label{lem5}
For positive integer $i,$ the characteristic polynomial of threshold graph $G^{\prime}=(0 1^{2i}0^{i+1} 1^{2i} 0^{2i+1} 1^{2i+2})$ is, to within a sign,
\begin{equation}
\label{eq8}
P_{G^{\prime}}(x) = x^{3i} (x+1)^{6i} (x+2i+2)(x^4 -(8i+2)x^3 -(-8i^2+4i+3)x^2 -(-8i^3 -20i^2-8i)x -8i^4-12i^3-4i^2)
 \end{equation}
\end{Lem}

\begin{Thr}
\label{main4}
The $n$-vertex graphs  $G=(0 1^{2i+1}0^i 1^{2i+2} 0^{2i+1} 1^{2i})$ and  $G^{\prime}=(0 1^{2i}0^{i+1} 1^{2i} 0^{2i+1} 1^{2i+2})$ are noncospectral and equienergetic graphs  if $n=9i+5, i\geq 1.$
\end{Thr}

\begin{Cor}
\label{cor2}
The $n$-vertex graphs  $G=(0 1^{2i+1}0^i 1^{2i+2} 0^{2i+1} 1^{2i})$ and  $G^{\prime}=(0 1^{2i}0^{i+1} 1^{2i} 0^{2i+1} 1^{2i+2})$ have the same energy but differ to the complete graph's energy. 
Furthermore, $E(G) = E(G^{\prime}) < E(K_n) = 18i +8,$  if $n=9i+5, i\geq 1.$
\end{Cor}

\end{document}